\documentclass[final,leqno]{article}
\usepackage{amsmath,amssymb}

\usepackage{graphicx}
\usepackage{dcolumn}
\usepackage{bm}
\usepackage{subeqnar}
\usepackage{epstopdf,overpic,color,rotating}
\usepackage{overpic}
\usepackage{gensymb}

\usepackage[paperheight=11in,paperwidth=8.5in]{geometry}

\graphicspath{{figures/}}

\newcommand*\samethanks[1][\value{footnote}]{\footnotemark[#1]}

\title{Multi-Resolution Dynamic Mode Decomposition}

\author{J. Nathan Kutz$^\ddagger$\thanks{Department of Applied Mathematics, University of Washington, Seattle, WA. 98195-2420.  $^\ddagger$ ({kutz@uw.edu}). Questions, comments, or corrections to this document may be directed to that email address.  A video summary/abstract of this work may be found at:  http://youtu.be/E1dNE02LaCE}
 \and Xing Fu\samethanks[1]
\and Steven L. Brunton\thanks{Department of Mechanical Engineering, University of Washington, Seattle, WA. 98195-2420.}
}

\begin{document}
\maketitle

\begin{abstract}
We demonstrate that the integration of the recently developed dynamic mode decomposition (DMD) with
a multi-resolution analysis allows for a decomposition method capable of robustly separating complex systems into a hierarchy of multi-resolution time-scale components.  A one-level separation allows for background (low-rank) and foreground (sparse) separation of dynamical data, or robust principal component analysis.  The multi-resolution dynamic mode decomposition is capable of  characterizing nonlinear dynamical systems in an equation-free manner by recursively decomposing the state of the system into low-rank terms whose temporal coefficients in time are known.  DMD modes with temporal frequencies near the origin (zero-modes) are interpreted as background (low-rank) portions of the given dynamics, and the terms with temporal frequencies bounded away from the origin are their sparse counterparts.  The multi-resolution dynamic mode decomposition (mrDMD) method is demonstrated on several examples involving multi-scale dynamical data, showing
excellent decomposition results, including sifting the El Ni\~no mode from ocean temperature data.
It is further applied to decompose a video data set into separate objects moving at different rates against a slowly varying background.  These examples show that the decomposition is an effective dynamical systems tool for data-driven discovery.
\end{abstract}



\section{Introduction}

Modeling of multi-scale systems, both in time and space, pervade modern theoretical and computational efforts across the engineering, biological and physical sciences.  Driving innovations are methods and algorithms that circumvent the significant challenges in efficiently connecting micro-scale to macro-scale effects that are separated potentially by orders of magnitude spatially and/or temporally.   Wavelet-based methods and/or windowed Fourier Transforms are ideally structured to perform such multi-resolution analyses (MRA) as they systematically remove temporal or spatial features by a process of recursive refinement of sampling from the data of interest~\cite{wavelet,wavelet2,DataBook}.  Typically, MRA is performed in either space or time, but not both simultaneously.  We propose integrating the concept of MRA with the recently developed Dynamic Mode Decomposition (DMD)~\cite{DMD0,DMD1,DMD4,DMD2,DMD3,DMD5}, a technique that produces low-dimensional spatio-temporal modes.  The proposed multi-resolution DMD (mrDMD) is shown to naturally integrate space and time so that the multi-scale spatio-temporal features are readily separated and
approximate dynamical models constructed.

The origins of the DMD method, which arose from pioneering work connecting the Koopman operator to 
dynamical systems theory~\cite{Mezic2004,Mezic2005}, are associated with the fluid dynamics community and the modeling of complex flows~\cite{DMD1,DMD4}.  Its growing success stems from the fact that it is an {\em equation-free}, data-driven method~\cite{DataBook} capable of providing accurate assessments of the spatio-temporal coherent structures in a complex system, or short-time future estimates, thus potentially allowing for control protocols to be enacted simply from data sampling. 
The mathematical architecture advocated here is an alternative to the equation-free, multi-scale modeling method
proposed by Kevrekidis and co-workers~\cite{kev1,kev2}.
More broadly, DMD has quickly gained popularity since it provides information about nonlinear dynamical systems.  DMD analysis can be considered to be a numerical approximation to Koopman spectral analysis~\cite{DMD4,mezic2}, and it is in this sense that DMD is applicable to nonlinear systems. In fact, the terms {\em DMD mode} and {\em Koopman mode} are often used interchangeably in the fluids literature.

At its core, the DMD method can be thought of as an ideal combination of the Proper Orthogonal Decomposition (POD), a spatial dimensionality-reduction technique, with Fourier Transforms in time.  
More precisely, the DMD method produces a least-square regression to
a linear dynamical system over the range of data collection.
In this work, our goal is to integrate the DMD decomposition with key concepts from wavelet theory and MRA.   Specifically, the DMD method takes snapshots of an underlying dynamical system to construct its decomposition.  However, the frequency and duration (sampling window) of the data collection process can be adapted, much as in wavelet theory, to sift out information at different scales.  Indeed, an iterative refinement of progressively shorter snapshot sampling windows and recursive extraction of DMD modes from slow to increasingly fast time scales allows for the mrDMD.  Moreover, it also allows for improved analytic predictions of the short-time future state of the system which is of critical importance for feedback control, for instance.  Critical innovations demonstrated here are the ability of the mrDMD to handle transient phenomenon and moving (translating/rotating) structures in data, both weaknesses of SVD-based techniques.

The paper is outlined as follows:  In Sec.~\ref{sec:dmd} the basic DMD theory is outlined with an emphasis on its low-rank approximation of data.  This is followed in Sec.~\ref{sec:mrDMD} by the development of the mrDMD structure and algorithm used in the subsequent applications, described in Sec.~\ref{sec:app}.  The applications exhibit the application to dynamical systems as well as more broadly to other applications such as video analysis.  The paper is concluded in Sec.~\ref{sec:discussion} with an overview and outlook of the method.

\section{Dynamic Mode Decomposition}
\label{sec:dmd}

The DMD method provides a spatio-temporal decomposition of data into
a set of dynamic modes that are derived from snapshots or measurements of a given system in time.   The mathematics underlying the extraction of dynamic information from time-resolved snapshots  is closely related to the idea of the Arnoldi algorithm~\cite{DMD1}, one of the workhorses of fast computational solvers.   The data collection process involves two (integer) parameters:
\begin{subeqnarray}
  &&  N = \mbox{number of spatial measurements per time snapshot, } \\
  &&  M= \mbox{number of snapshots taken in time.}
\end{subeqnarray}
Originally the algorithm was designed to collect data at regularly spaced intervals of time.  However, new innovations, and a more general definition of the DMD, allow for both sparse spatial~\cite{cdmd} and temporal~\cite{Tu2014ef} collection of data as well as irregularly spaced collection times~\cite{DMD5}.  To illustrate the algorithm, we consider regularly spaced sampling in time.  The data collection times are given by:
\begin{equation}
t_{m+1} = t_{m} + \Delta t 
\end{equation}
where the collection time starts at $t_1$ and ends at $t_M$, and the
interval between data collection times is $\Delta t$.  In the mrDMD method, the total number of snapshots will vary as the algorithm extracts multi-timescale spatio-temporal structures.  This will be the central focus of the next section.

The data snapshots are arranged into an $N\times M$ matrix
\begin{equation}
  {\bf X} = \left[ {\bf x}(t_1) \,\,\, {\bf x}(t_2) \,\,\, {\bf x}(t_3) \,\,\, \cdots \,\,\,  {\bf x}(t_M)    \right]
\end{equation}
where the vector ${\bf x}$ are the $N$ measurements of the state variable of the system of interest at the data collection points.   The objective is to mine the data matrix ${\bf X}$ for important dynamical information.  For the purposes of the DMD method, the following matrix is also defined:
\begin{equation}
  {\bf X}_j^{k} = \left[  {\bf x}(t_j) \,\,\, {\bf x}(t_{j+1}) \,\,\,  \cdots \,\,\, {\bf x}(t_k)    \right].
\end{equation}
Thus this matrix includes columns $j$ through $k$ of the original data matrix.

The DMD method approximates the modes of the so-called {\em Koopman operator}.
The Koopman operator is a linear, infinite-dimensional operator that represents
nonlinear, possibly infinite-dimensional, dynamics without linearization~\cite{DMD4,mezic2},
and it is the adjoint of the Perron-Frobenius operator\index{Perron-Frobenius operator}.  
The method can be viewed as computing, from the experimental data, the eigenvalues and eigenvectors (low-dimensional
modes) of a linear model that approximates the underlying dynamics, even if the dynamics
are nonlinear.  Since the model is assumed to be linear, the decomposition gives
the growth rates and frequencies associated with each mode.  If the underlying model
is linear, then the DMD method recovers the leading eigenvalues and eigenvectors 
normally computed using standard solution methods for linear differential equations.

The DMD involves approximating the eigendecomposition of the best-fit linear operator ${\bf A}$ that relates a state ${\bf x}_j$ at time $t_j$ to the state ${\bf x}_{j+1}$ at the next timestep: 
\begin{equation}
   {\bf x}_{j+1} \approx {\bf A} {\bf x}_j.
   \label{eq:Koopman}
\end{equation}
If Eq.~\eqref{eq:Koopman} holds exactly and the data is generated by a linear system, then each column of ${\bf X}_1^{M-1}$ is an element of a Krylov subspace:
\begin{equation}
  {\bf X}_{1}^{M-1} = \left[ {\bf x}_1 \,\,\, {\bf A}{\bf x}_1 \,\,\, {\bf A}^2{\bf x}_1 \,\,\, \cdots \,\,\, 
    {\bf A}^{M-2}{\bf x}_{1} \right] \, .
    \label{eq:kry1}
\end{equation}
Eq.~\eqref{eq:Koopman} may be written in matrix form as
\begin{equation}
   {\bf X}_{2}^M \approx {\bf A} {\bf X}_1^{M-1},
   \label{eq:Koopman2}
\end{equation}
where the operator ${\bf A}$ is chosen to minimize the Frobenius norm of $\|{\bf X}_2^M-{\bf A}{\bf X}_1^{M-1}\|_F$.  In other words, the operator ${\bf A}$ advances each snapshot column in ${\bf X}_1^{M-1}$ a single timestep, $\Delta t$, resulting in the future snapshot columns in ${\bf X}_2^M$.

In practice, when the state dimension $N$ is large, the matrix ${\bf A}$ may be intractable to analyze directly.  Instead, DMD circumvents the eigendecomposition of ${\bf A}$ by considering a rank-reduced representation in terms of a POD-projected matrix $\tilde{\bf A}$.  The DMD algorithm proceeds as follows~\cite{DMD5}:

\begin{enumerate}
\item First, take the SVD of ${\bf X}_1^{M-1}$\cite{tre}:
\begin{equation}
  {\bf X}_{1}^{M-1} = {\bf U} {\bf \Sigma} {\bf V}^*,
  \label{eq:kry3}
\end{equation}
where $*$ denotes the conjugate transpose, 
${\bf U}\in {\mathbb{C}}^{N\times K}$, ${\bf \Sigma}\in {\mathbb{C}}^{K\times K}$
and ${\bf V}\in {\mathbb{C}}^{M-1\times K}$.  Here $K$ is the rank of the reduced SVD
approximation to ${\bf X}_{1}^{M-1}$.  The left singular vectors ${\bf U}$ are POD modes.

The SVD reduction in (\ref{eq:kry3}) could also be exploited at this stage in the algorithm 
to perform a low-rank truncation of the data.  Specifically, if low-dimensional structure is
present in the data, the singular values of ${\bf \Sigma}$ will decrease sharply to zero with
perhaps only a limited number of dominant modes.   A principled way to truncate noisy data would be to use the
recent hard-thresholding algorithm of Gavish and Donoho~\cite{gavish}.

\item Next, compute $\tilde{\bf A}$, the $K\times K$ projection of the full matrix ${\bf A}$ onto POD modes:
\begin{eqnarray}
{\bf A} &=& {\bf X}_2^M{\bf V}\boldsymbol{\Sigma}^{-1}{\bf U}^*\nonumber\\
\Longrightarrow\quad\tilde{\bf A} &=& {\bf U}^*{\bf A}{\bf U} = {\bf U}^*{\bf X}_2^M{\bf V}\boldsymbol{\Sigma}^{-1}.
\end{eqnarray}
\item Compute the eigendecomposition of $\tilde{\bf A}$:
\begin{equation}
\tilde{\bf A}{\bf W} = {\bf W}\boldsymbol{\Lambda},
\end{equation}
where columns of ${\bf W}$ are eigenvectors and $\boldsymbol{\Lambda}$ is a diagonal matrix containing the corresponding eigenvalues $\lambda_k$.
\item Finally, we may reconstruct eigendecomposition of $\bf A$ from ${\bf W}$ and $\boldsymbol{\Lambda}$.  In particular, the eigenvalues of $\bf A$ are given by $\boldsymbol{\Lambda}$ and the eigenvectors of ${\bf A}$ (DMD modes) are given by columns of $\bf \Psi$:
\begin{equation}
{\bf \Psi} = {\bf X}_2^M{\bf V}\boldsymbol{\Sigma}^{-1}{\bf W}.
\label{eq:DMDjtu}
\end{equation}
\end{enumerate}
Note that Eq.~\eqref{eq:DMDjtu} from~\cite{DMD5} differs from the formula ${\bf \Psi}={\bf U}{\bf W}$ from~\cite{DMD1}, although these will tend to converge if ${\bf X}_1^{M-1}$ and ${\bf X}_2^M$ have the same column spaces.
With the low-rank approximations of both the eigenvalues and eigenvectors in
hand, the projected future solution can be constructed for all time in the future.
By first rewriting for convenience $\omega_k=\ln(\lambda_k)/\Delta t$ , then the approximate solution
at all future times, $ {\bf x}_{\mbox{\tiny DMD}}(t)$, is given by
\begin{equation}
  {\bf x}_{\mbox{\tiny DMD}}(t) = \sum_{k=1}^{K} b_k(0) \psi_k ({\boldsymbol{\xi}}) \exp(\omega_k t) = {\bf \Psi} 
  \mbox{diag} (\exp(\omega t)) {\bf b} 
  \label{eq:dmd_sol}
\end{equation}
where ${\boldsymbol{\xi}}$ are the spatial coordinates, 
$b_k(0)$ is the initial amplitude of each mode, ${\bf \Psi}$ is the matrix whose
columns are the eigenvectors $\psi_k$, $\mbox{diag}(\omega t)$ is a diagonal
matrix whose entries are the eigenvalues $\exp(\omega_k t)$, and ${\bf b}$ is a
vector of the coefficients $b_k$.\index{eigenvalues}\index{eigenfunctions}\index{POD}

An alternative interpretation of (\ref{eq:dmd_sol}) is that it represents the least-square
fit, or regression, of a linear dynamical system $d{\bf x}_{\mbox{\tiny DMD}}/dt = {\bf A} {\bf x}_{\mbox{\tiny DMD}}$
to the data sampled.  In particular, the matrix ${\bf A}$ constructed is such that
$\|{\bf x}(t)-{\bf x}_{\mbox{\tiny DMD}}(t)\|$ is minimized.  In the context of the multi-resolution analysis
that follows, each level of the multi-scale decomposition produces a linear dynamical system, or matrix ${\bf A}$, for the
time-scale under consideration.

It only remains to compute the initial coefficient values $b_k(0)$.  If we consider the
initial snapshot (${\bf x}_1$) at time $t_1=0$, let's say, then (\ref{eq:dmd_sol}) gives
${\bf x}_1 ={\bf \Psi} {\bf b}$.  This generically is not a square matrix so that
its solution 
\begin{equation}
  {\bf b} = {\bf \Psi}^+ {\bf x}_1
\end{equation}
can be found using a pseudo-inverse.  Indeed, ${\bf \Psi}^+$ denotes the Moore-Penrose
pseudo-inverse\index{pseudo-inverse} that can be accessed in MATLAB via the \texttt{pinv} command.  The pseudo-inverse is equivalent to
finding the best solution ${\bf b}$ the in the least-squares (best fit) sense.  This is equivalent
to how DMD modes were derived originally.

Overall then, the DMD algorithm presented here takes advantage of low dimensionality
in the data in order to make a low-rank approximation of the linear mapping that
best approximates the nonlinear dynamics of the data collected for the system.  Once
this is done, a prediction of the future state of the system is achieved for all time.  Unlike
the POD-Galerkin method, which requires solving a low-rank set of dynamical quantities to
predict the future state, no additional work is required for the future state prediction outside
of plugging in the desired future time into (\ref{eq:dmd_sol}).  Thus the advantages
of DMD revolve around the fact that (i) it is an equation-free architecture, and (ii) a future state
prediction is known for any time $t$ (of course, provided the DMD approximation holds).

More broadly, the DMD method was shown to be a highly successful method for
foreground/background subtraction in video feeds~\cite{jake}.  Indeed, the DMD method
is a novel, dynamical systems base method for performing a  {\em robust Principal Components Analysis} 
(RPCA) of data streams~\cite{jake}.   Importantly, the DMD-based RPCA performs
the low-rank/sparse matrix separation 3-4 orders of magnitude faster than standard
$\ell_1$ optimization methods~\cite{candes}.  RPCA is extremely important for handling
data outliers and/or corrupt data matrices.  In the multi-resolution version of DMD, the
RPCA can be used to effectively remove outliers at each level of decomposition.

%
%
%
%
%
%

\section{Multi-Resolution Dynamic Mode Decomposition}
\label{sec:mrDMD}

The mrDMD is inspired by the observation that the slow- and fast-modes can be separated for such applications as foreground/background subtraction in video feeds~\cite{jake}.  The mrDMD recursively removes low-frequency, or slowly-varying, content from a given collection of snapshots.  Typically, the number of snapshots $M$ are chosen so that the DMD modes provide an {approximately full rank approximation} of the dynamics observed.  Thus $M$ is chosen so that all high- and low-frequency content is present.  In the mrDMD, $M$ is originally chosen in the same way so that an {approximate full rank approximation} can be accomplished.  However, from this initial pass through the data, the slowest $m_1$ modes are removed, and the domain is divided into two segments with $M/2$ snapshots each.  DMD is once again performed on each $M/2$ snapshot sequences.  Again the slowest $m_2$ modes are removed and the algorithm is continued until a desired termination.

Mathematically, the mrDMD separates the DMD approximate solution (\ref{eq:dmd_sol}) in the first pass as follows:
\begin{eqnarray}\hspace*{.3in}
{\bf x}_{\mbox{\tiny mrDMD}}(t) \!=\!\!\! \sum_{k=1}^{M} \! b_k(0) \psi_k^{(1)} \!({\boldsymbol{\xi}}) \exp(\omega_k t) \! &=& \!\!\! \sum_{k=1}^{m_1} \! b_k(0) \psi_k^{(1)}  \! ({\boldsymbol{\xi}}) \exp(\omega_k t)
  \! + \!\!\!\!\!\!\! \sum_{k=m_1+1}^{M} \!\! b_k(0) \psi_k^{(1)} \! ({\boldsymbol{\xi}}) \exp(\omega_k t)
  \label{eq:dmd1} \\
  & & \,\,\, \mbox{(slow modes)} \hspace*{1.0in} \mbox{(fast modes)}  \nonumber
\end{eqnarray}
where the $\psi_k^{(1)}\! ({\bf x})$ represent the DMD  modes computed from the full $M$ snapshots.  

The first sum in this expression (\ref{eq:dmd1}) represents the slow-mode dynamics whereas the second sum is everything else.  Thus the second sum can be computed to yield the matrix:
\begin{equation}
   {\bf X}_{M/2} =  \sum_{k=m_1+1}^{M} \!\! b_k(0) \psi_k^{(1)} \! ({\boldsymbol{\xi}}) \exp(\omega_k t) \, .
\end{equation}
The DMD analysis outlined in the previous section can now be performed once again on the data matrix ${\bf X}_{M/2}$.  However, the matrix ${\bf X}_{M/2}$ is now separated into two matrices
\begin{equation}
   {\bf X}_{M/2} = {\bf X}_{M/2}^{(1)} + {\bf X}_{M/2}^{(2)} 
\end{equation}
where the first matrix contains the first $M/2$ snapshots and the second matrix contains the remaining $M/2$ snapshots.  The $m_2$ slow-DMD modes at this level are given by $\psi_k^{(2)}$, where they are computed separately in the first of second interval of snapshots.   

The iteration process works by recursively removing slow frequency components and building the new matrices ${\bf X}_{M/2}, {\bf X}_{M/4}, {\bf X}_{M/8}, \cdots$ until a desired/prescribed multi-resolution decomposition has been achieved.  The approximate DMD solution 
can then be constructed as follows:
\begin{eqnarray} \hspace*{.3in}
  {\bf x}_{\mbox{\tiny mrDMD}}(t) \! &=& \!   \sum_{k=1}^{m_1} \! b_k^{(1)} \psi_k^{(1)}  \! ({\boldsymbol{\xi}}) \exp(\omega_k^{(1)} t)
   + \! \sum_{k=1}^{m_2} \! b_k^{(2)} \psi_k^{(2)}  \! ({\boldsymbol{\xi}}) \exp(\omega_k^{(2)} t) \nonumber \\
   && + \! \sum_{k=1}^{m_3} \! b_k^{(3)} \psi_k^{(3)}  \! ({\boldsymbol{\xi}}) \exp(\omega_k^{(3)} t) + \cdots
  \label{eq:dmd2} 
\end{eqnarray}
where at the evaluation time $t$, the correct modes from the sampling window are selected at each level of the decomposition.
Specifically, the  $\psi_k^{(k)}$ and $\omega_k^{(k)}$ are the DMD modes and DMD eigenvalues at the $k$th level of decomposition, the $b_k^{(k)}$ are the initial projections of the data onto the time interval of interest, and the $m_k$ are the number of slow-modes retained at each level.  The advantage of this method is readily apparent:  different spatio-temporal DMD modes are used to represent key multi-resolution features.  Thus there is not a single set of modes that dominates the SVD and potentially marginalizes features at other time scales. 

\begin{figure}[t]
\vspace*{-.8in}
\begin{center}
\begin{overpic}[width=0.7\textwidth]{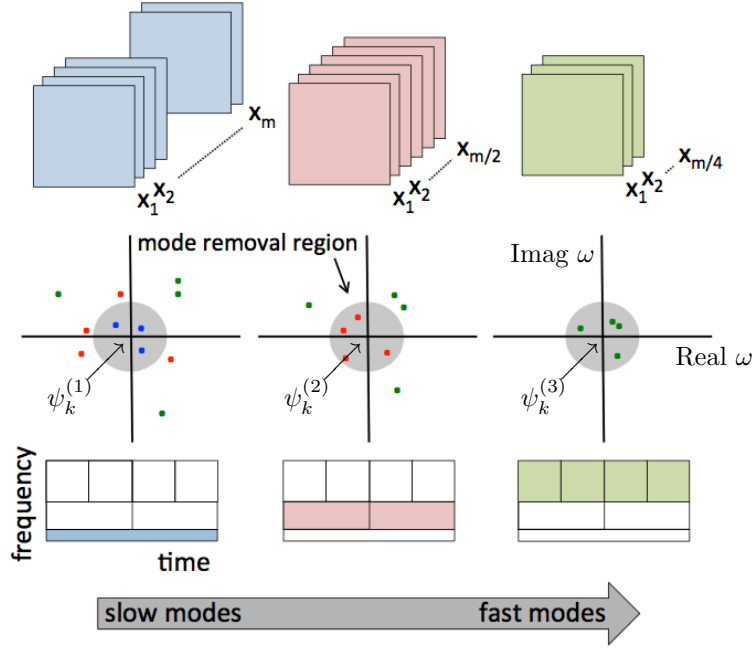}
\normalsize
\put(70,46){Real $\omega$}
\put(54,56){Imag $\omega$}
\put(55,42){$\psi^{(3)}_k$}
\put(58,44){\rotatebox{45}{$\xrightarrow{\hspace*{0.5cm}}$}}
\put(32,42){$\psi^{(2)}_k$}
\put(35,44){\rotatebox{45}{$\xrightarrow{\hspace*{0.5cm}}$}}
\put(9,42){$\psi^{(1)}_k$}
\put(12,44){\rotatebox{45}{$\xrightarrow{\hspace*{0.5cm}}$}}
\end{overpic}
\end{center}
\vspace*{-1.2in}
\caption{\label{fig:mrDMD} Representation of the multi-resolution dynamic mode decomposition where successive sampling of the data, initially with $M$ snapshots and decreasing by a factor of two at each resolution level, is shown (top figures).  The DMD spectrum is shown in the middle panel where there are $m_1$ (blue dots) slow-dynamic modes at the slowest level, $m_2$ (red) modes at the next level and $m_3$ (green) modes at the fastest time-scale shown.  The shaded region represents the modes that are removed at that level.  The bottom panels shows the wavelet-like time-frequency decomposition of the data color coded with the snapshots and DMD spectral representations.}
\end{figure}

Figure~\ref{fig:mrDMD} illustrates the multi-resolution DMD process pictorially.  In the figure, a three-level decomposition is performed with  the slowest scale represented in blue (eigenvalues and snapshots), 
the mid-scale in red and the fast scale in green.   The connection to multi-resolution wavelet analysis is also evident from the bottom panels as one can see that the mrDMD method successively pulls out time-frequency information in a principled way.  

As a final remark, the sampling strategy and algorithm discussed here (See Fig.~\ref{fig:mrDMD}) can be easily 
modified since only the slow modes at
each decomposition level need to be accurately computed.  Thus one can modify the algorithm so
as to sample a fixed number, for instance $M$, data snapshots in each sampling window.  The value of
$M$ need not be large as only the slow modes need to be resolved.  Thus
the sampling rate would increase as the decomposition proceeds from one level to the next.
This assures that the lowest levels of the mrDMD are not highly-sampled since the cost of the SVD
would be greatly increased by such a fine sampling rate.   
 
\subsection{Formal mrDMD Expansion}

The solution (\ref{eq:dmd2}) can be made more precise.  Specifically, one must account
for the number of levels ($L$) of the decomposition, the number of time bins ($J$) for each level,
and the number of modes retained at each level ($m_L$).  This can be easily seen
in Fig.~\ref{fig:mrDMD}.  Thus the solution is parametrized
by the following three indices:
\begin{subeqnarray}
  &&  \ell = 1, 2, \cdots , L \,\,\,  \mbox{number of decomposition levels } \\
  &&  j = 1, 2, \cdots, J  \,\,\, \mbox{number time bins per level} \,\, (J=2^{(\ell-1)})  \\
  &&  k= 1, 2, \cdots, m_{L} \,\,\, \mbox{number of modes extracted at level $L$.}
\end{subeqnarray}
To formally define the series solution for ${\bf x}_{\mbox{\tiny mrDMD}}(t)$, we
propose the following indicator function
\begin{equation}
  f_{\ell,j} (t) = \left\{  \begin{array}{cl}  1 & t\in[t_j,t_{j+1}] \\ 0 & \mbox{elsewhere}  \end{array}    \right.
  \,\,\, \mbox{with} \,\,\,\,\, j=1,2,\cdots,J \,\,\,\,\, \mbox{and} \,\,\,\,\, J=2^{(\ell-1)}
  \label{eq:indicator}
\end{equation}
which is only non-zero in the interval, or time bin, associated with the value of $j$.  The parameter
$\ell$ denotes the level of the decomposition.

\begin{figure}[t]
\begin{center}
\vspace*{-2.5in}
\begin{overpic}[width=1.0\textwidth]{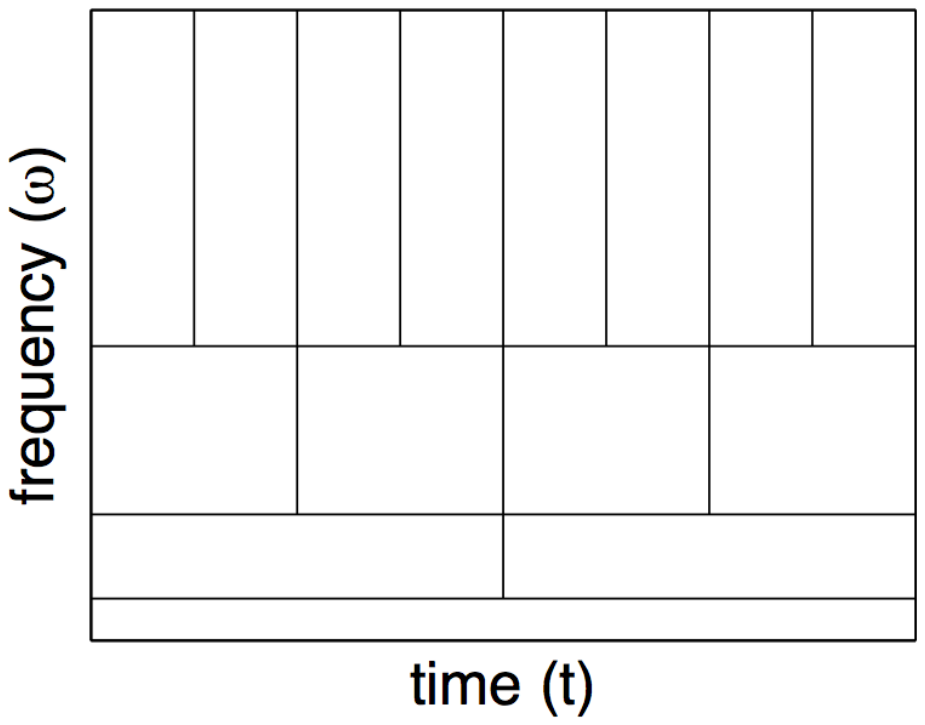}
\put(50,50){\Large{$\psi_k^{(\ell,j)}({\boldsymbol{\xi}})$}} 
\put(51.4,45){\rotatebox{90}{$\xrightarrow{\hspace*{0.5cm}}$}}
\put(50,43){$k=$ mode number at level ${\ell}$}
\put(49.8,56){\rotatebox{-60}{$\xrightarrow{\hspace*{0.5cm}}$}}
\put(47,58){$\ell=$ decomposition level}
\put(54,55.5){\rotatebox{-160}{$\xrightarrow{\hspace*{1.0cm}}$}}
\put(61,54.5){$j=$ time bin}
\put(40,40){$\psi_k^{(1,1)}$}
\put(11.7,42.2){$\psi_k^{(2,1)}$}
\put(27,42.2){$\psi_k^{(2,2)}$}
\put(8,47){$\psi_k^{(3,1)}$}
\put(15.5,47){$\psi_k^{(3,2)}$}
\put(23,47){$\psi_k^{(3,3)}$}
\put(30.5,47){$\psi_k^{(3,4)}$}
\put(6,65){$\psi_k^{(4,1)} \,\,\, \cdots$}
\put(32,65){$\psi_k^{(4,8)}$}
\put(7,64){\rotatebox{-90}{$\xrightarrow{\hspace*{0.8cm}}$}}
\put(33,64){\rotatebox{-90}{$\xrightarrow{\hspace*{0.8cm}}$}}
\put(33,41){\rotatebox{-180}{$\xrightarrow{\hspace*{1.0cm}}$}}
\normalsize
\end{overpic}
\vspace*{-3in}
\end{center}
\caption{\label{fig:wavelet} Illustration of the mrDMD mode decomposition and
hierarchy.  Represented are the modes $\psi_k^{\ell,j}({\boldsymbol{\xi}})$ and their position in the decomposition
structure.  The triplet of integer values, $\ell, j$ and $k$, uniquely express the time level, bin and mode of
the decomposition.}
\end{figure}

The three indices and  indicator function (\ref{eq:indicator}) give the mrDMD solution expansion
\begin{equation}
  {\bf x}_{\mbox{\tiny mrDMD}}(t)  = \sum_{\ell=1}^{L} \sum_{j=1}^{J} \sum_{k=1}^{m_L}
  f_{\ell,j}(t)  b_k^{(\ell,j)} \psi_k^{(\ell,j)}  \! ({\boldsymbol{\xi}}) \exp(\omega_k^{(\ell,j)} t) 
  \, .
  \label{eq:dmd3}
\end{equation}
This is a concise definition of the mrDMD solution that includes the information on the level,
time bin location and number of modes extracted.
Figure~\ref{fig:wavelet} demonstrates the mrDMD decomposition in terms of the solution
(\ref{eq:dmd3}).  In particular, each mode is represented in its respective time bin and level.
An alternative interpretation of this solution is that it yields the least-square fit, at each
level $\ell$ of the decomposition, to the linear dynamical system
\begin{equation}
  \frac{d {\bf x}^{(\ell,j)}}{dt} = {\bf A}^{{(\ell,j)}} {\bf x}^{(\ell,j)} \, 
\end{equation}
where the matrix ${\bf A}^{(\ell,j)}$ captures the dynamics in a given time bin $j$ at level ${\ell}$.

The indicator function $f_{\ell,j}(t)$ acts as sifting function for each time bin.  Interestingly,
this function acts as the Gab\'or window of a windowed Fourier transform~\cite{DataBook}.  
Since our sampling bin has a hard cut-off of the time series, it may introduce some artificial
high-frequency oscillations.  Time-series analysis, and wavelets in particular, introduce
various functional forms that can be used in an advantageous way.  Thus thinking 
more broadly, one can imagine using wavelet functions for the sifting operation, thus allowing
the time function $f_{\ell,j}(t)$ to take the form of one of the many potential wavelet basis, i.e.
Haar, Daubechies, Mexican Hat, etc.  This will be considered in future work.  For the present,
we simply use the sifting function introduced in (\ref{eq:indicator})

\section{Application of Method}
\label{sec:app}

The  mrDMD developed in the last section is implemented here on three 
example data sets:  the first being an exemplar of a video stream, the second
coming from atmospheric-ocean data, and the third an exemplar of data with moving objects.  This
final example is especially important as it renders many SVD-based decompositions useless.  
For the first case, the implementation is compared against
the standard DMD algorithm, highlighting the ability of the mrDMD to correctly 
capture multi-scale phenomena.   It should be noted that in previous work with DMD, which was not framed in
the mrDMD architecture advocated here, a level-1 decomposition was
effectively used to separate the foreground from the background in a video~\cite{jake}.  Indeed, this work motivates the generalization to the
multi-resolution analysis.

\begin{figure}[t]
\begin{center}
\begin{overpic}[width=0.8\textwidth]{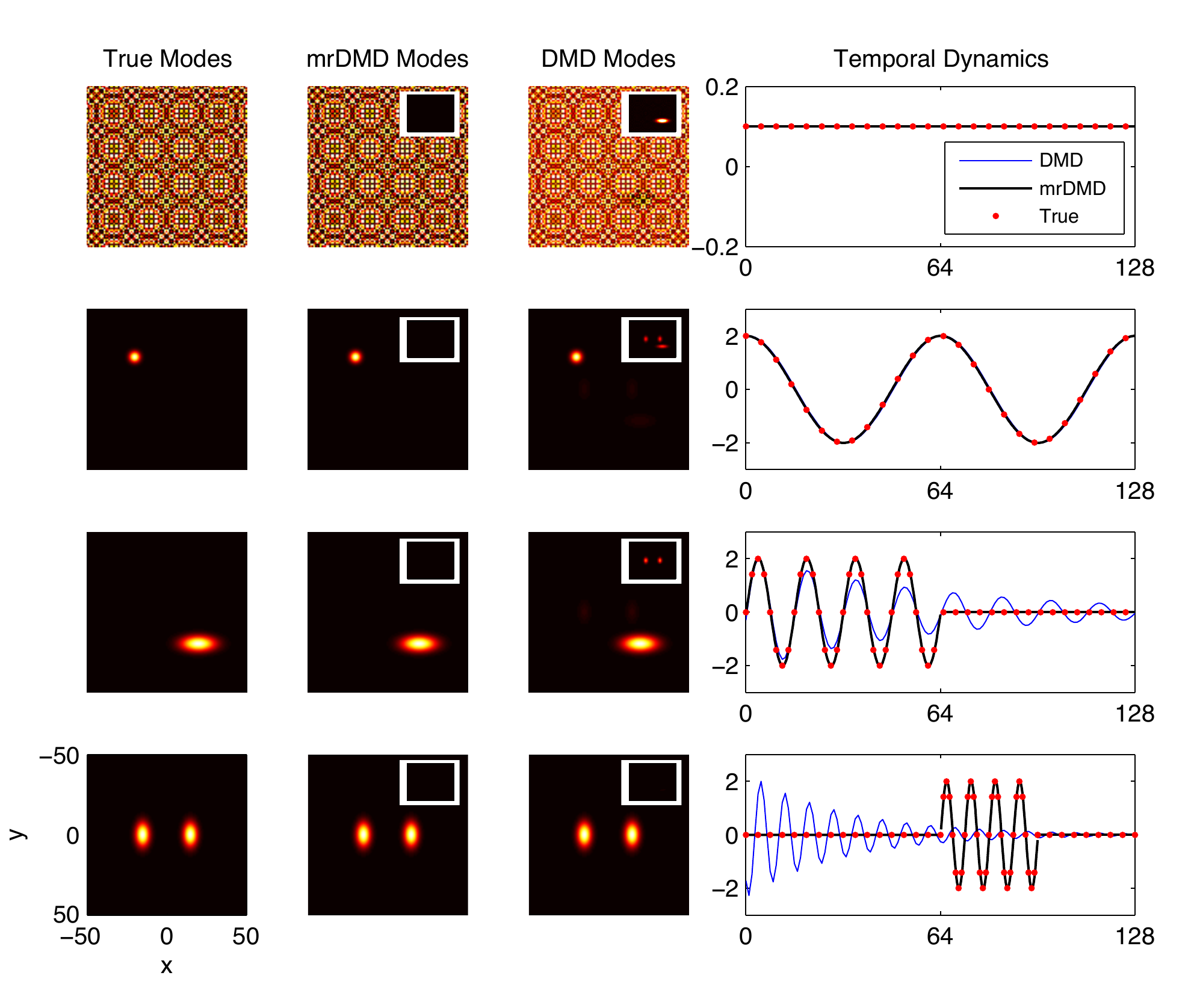}
\put(3,72){$\bar{\psi}_1$}
\put(3,53.5){$\bar{\psi}_2$}
\put(3,35.5){$\bar{\psi}_3$}
\put(3,17){$\bar{\psi}_4$}
\put(30,60.2){${\psi}_1^{(1,1)}$}
\put(30,42){${\psi}_1^{(2,1)}$}
\put(30,23.3){${\psi}_1^{(3,1)}$}
\put(30,4.5){${\psi}_1^{(3,3)}$}
\put(50,60.2){${\psi}_1$}
\put(50,42){${\psi}_2$}
\put(50,23.3){${\psi}_3$}
\put(50,4.5){${\psi}_4$}
\put(65,75){$a_1(t)$}
\put(65,56.3){$a_2(t)$}
\put(65,37.7){$a_3(t)$}
\put(65,19){$a_4(t)$}
\put(49,18){\color{white}{$R_3$}}
\put(29,18){\color{white}{$R_{3,3}$}}
\put(49,36.5){\color{white}{$R_3$}}
\put(29,36.5){\color{white}{$R_{3,1}$}}
\put(49,56){\color{white}{$R_2$}}
\put(29,56){\color{white}{$R_{2,1}$}}
\put(49,74.5){\color{white}{$R_1$}}
\put(29,74.5){\color{white}{$R_{1,1}$}}
\normalsize
\end{overpic}
\end{center}
\caption{\label{fig:modes_temporal} Comparison of the true modes ($\bar{\psi}_k$) to DMD ($\psi_k$) and 
mrDMD ($\psi_k^{j}$) modes.  The true, mrDMD and DMD modes are showed in columns from left to right.  The inset shows the relative error of the mrDMD ($R_{j,\ell}=\| \bar{\psi}_j - \psi^{(\ell,1)}_j \|$) and DMD (${R}_{j}=\| \bar{\psi}_j - \psi_j \|$) in comparison to the true mode. Temporal dynamics are also compared, with true dynamics shown in red dots, DMD reconstructed dynamics shown in blue solid line and mrDMD reconstructed dynamics shown in black solid line.}
\end{figure}

\subsection{Spatio-temporal filtering of video}

\begin{figure}[t]
\begin{center}
\begin{overpic}[width=0.8\textwidth]{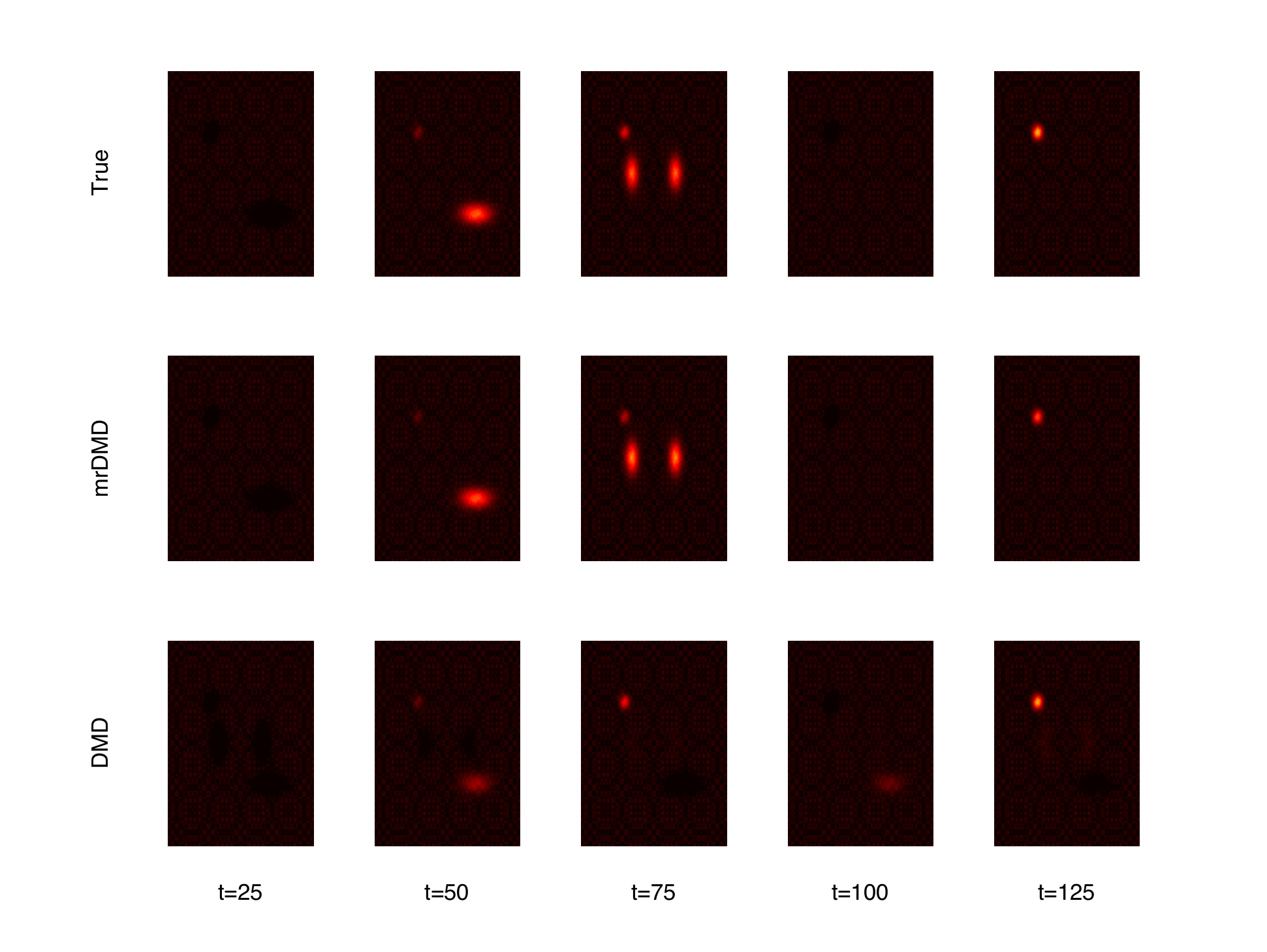}
\put(23,6){50}
\put(10,23){50}
\put(11,6){-50}
\put(9,8){-50}
\normalsize
\end{overpic}
\end{center}
\caption{\label{fig:compare} Comparison of snapshots from the original video of Fig.~\ref{fig:modes_temporal} and reconstruction from DMD and mrDMD.  The true dynamics are showed in the first row, while the mrDMD reconstructed dynamics are showed in the second row and the DMD reconstructed dynamics are showed in the third row.  The time labels on the bottom correspond to when each snapshot was taken.  The mrDMD reconstruction matches with the original video very well while the DMD reconstruction shows inconsistencies.  Indeed, the mrDMD method can easily handle signals that
turn on and off in time while the DMD method cannot.}
\end{figure}
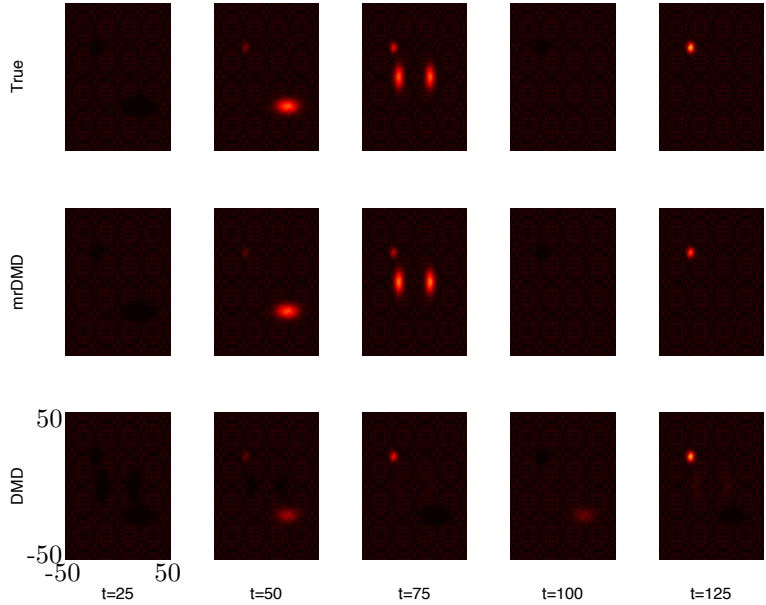

The first example we consider for application of the mrDMD is illustrated in Fig.~\ref{fig:modes_temporal}.
For this example, four spatio-temporal modes are combined into a single data set.   Specifically, we combine
the four modes shown in the left panels with the time dynamics given in the right panels.  The 
four modes used to construct the true solution are represented by $\bar{\psi}_j (x,y)$ for
$j=1,2,3$ and 4.  Their corresponding time dynamics are given by $a_j(t)$.  Thus the true solution
is expressed by 
\begin{equation}
  \bar{\bf x} = \sum_{j=1}^4  a_j(t) \bar{\psi}_j (x,y) \, .
  \label{eq:test}
\end{equation}
Both the DMD (represented by ${\bf x}_{\mbox{\tiny DMD}}$ and the modes $\psi_j$ of (\ref{eq:dmd_sol}) with $j=1, 2, 3$ and 4) and mrDMD (represented by ${\bf x}_{\mbox{\tiny mrDMD}}$ and
the $\psi_{k}^{(\ell,j)}$ of (\ref{eq:dmd2}) where $k=1$ and $\ell=1, 2, 3$) attempt to 
reconstruct $\bar{\bf x}$.  

\begin{figure}[t]
\begin{center}
\begin{overpic}[width=0.8\textwidth]{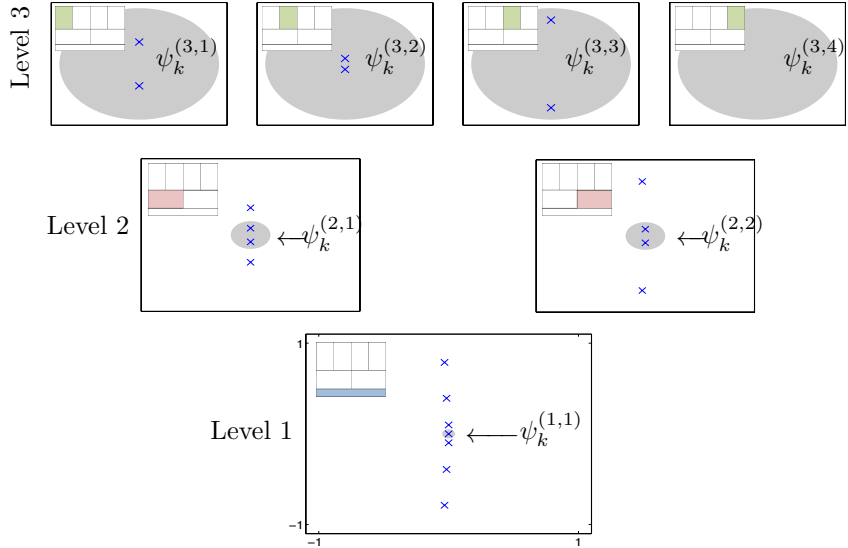}
\normalsize
\put(58,19.5){$\psi^{(1,1)}_k$}
\put(52,20.5){\rotatebox{180}{$\xrightarrow{\hspace*{0.5cm}}$}}
\put(34,41){$\psi^{(2,1)}_k$}
\put(31.2,42.2){\rotatebox{180}{$\xrightarrow{\hspace*{0.2cm}}$}}
\put(78,41){$\psi^{(2,2)}_k$}
\put(75.2,42.2){\rotatebox{180}{$\xrightarrow{\hspace*{0.2cm}}$}}
\put(18,60.5){$\psi^{(3,1)}_k$}
\put(41,60.5){$\psi^{(3,2)}_k$}
\put(63,60.5){$\psi^{(3,3)}_k$}
\put(87,60.5){$\psi^{(3,4)}_k$}
\put(24,19.5){Level 1}
\put(6,42){Level 2}
\put(2,58){\rotatebox{90}{Level 3}}
\end{overpic}
\end{center}
\vspace*{-.2in}
\caption{\label{fig:evals} Eigenvalues from level-1 to level-3 of the mrDMD decomposition applied to the data generated in Fig.~\ref{fig:modes_temporal}.  This is a specific case of the mrDMD abstraction shown in Fig.~\ref{fig:mrDMD}.  The shaded circle shows the threshold radius used for background mode subtraction.  The inset illustrates the time-frequency position of each sampling window.  Also illustrated are the mode selection and labeling process. }
\end{figure}

Figure~\ref{fig:modes_temporal} shows the true modes along with the approximating
mrDMD and DMD modes.  The inset in the mrDMD and DMD modes shows the $\ell_2$ difference
between the true modes and the approximating modes of the respective decomposition, i.e. $R_{j,\ell}=\| \bar{\psi}_j - \psi^{(\ell,1)}_j \|$ and ${R}_{j}=\| \bar{\psi}_j - \psi_j \|$ respectively.  The mrDMD modes are
almost identical to the true modes while the DMD modes show an error in mixing the modes.
The error in the DMD modes is clearly illustrated in the time dynamics.  Specifically, the DMD
does not correctly capture the on-off dynamics of modes three and four so that $\| \bar{\bf x} - {\bf x}_{\mbox{\tiny DMD}} \|  \sim O(1)$.  In contrast, the mrDMD
is able to easily capture the correct time dynamics with its windowing decomposition. Indeed,
the agreement between the mrDMD and exact solution is remarkable, i.e. $\| \bar{\bf x} - {\bf x}_{\mbox{\tiny mrDMD}} \| \ll 1$.

To further illustrate the accuracy of the mrDMD decomposition in comparison to the DMD decomposition,
consider Fig.~\ref{fig:compare}.  This shows the true solution $\bar{\bf x}$ at the time snap shots
of $t=25, 50, 75, 100$ and 125.  The mixing of the four modes of (\ref{eq:test}) is nicely captured
by the mrDMD while the DMD fails to capture key features at these time points.  Indeed, at time $t=75$
in particular, the agreement is quite poor between the DMD and true solution.

As a final part of the analysis of the video sequence represented by $\bar{\bf x}(t)$ in (\ref{eq:test}),
we can explicitly perform the multi-resolution decomposition depicted in Fig.~\ref{fig:mrDMD}.  In this case,
the dominant eigenvalue are kept  at each level of the decomposition.  The bottom centered box is the
level 1 decomposition showing that there exists a single mode in the mode removal region (shaded blue inset).
This is mode $\psi_1^{(1,1)}$.
At the next level, there are two modes (complex conjugate pairs) in the mode 
removal region (shaded pink inset), $\psi_1^{(2,1)}$ and $\psi_2^{(2,1)}$ (left panel) and 
$\psi_1^{(2,2)}$ and $\psi_2^{(2,2)}$ (right panel).  These modes represent the slow oscillatory dynamics driven by $a_2(t)$.  
And in the final level, an additional pair of eigenvalues are present in the first three frames of the decomposition
that correspond to the dynamics of modes $a_3(t)$ (frames 1 and 2) and  $a_4(t)$ (frame 3) with modes
$\psi_1^{(3,j)}$ and $\psi_2^{(3,j)}$.  Thus
the multi-resolution strategy outlined in Fig.~\ref{fig:mrDMD} holds.  Moreover, the intuition derived
from this wavelet-like strategy is remarkably effective in decomposing the multi-resolution, spatio-temporal
dynamics of $\bar {\bf x}$.

\subsection{Multi-scale time separation of complex system:  El Ni\~no, Southern Oscillation}

\begin{figure}[t]
\begin{center}
\begin{overpic}[width=0.8\textwidth]{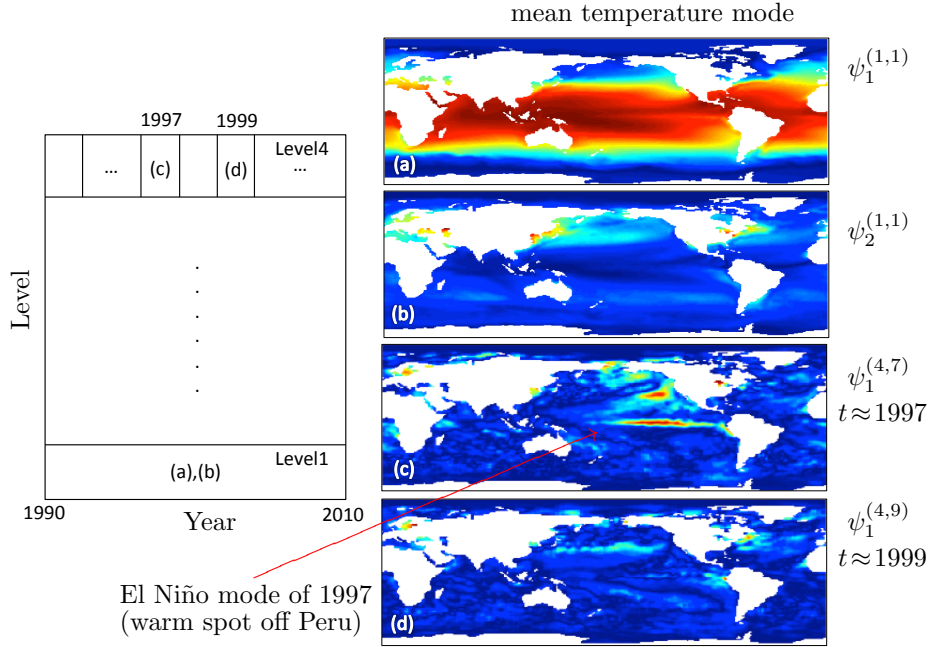}
\normalsize
\put(30,12){\color{red}{\rotatebox{23}{$\xrightarrow{\hspace*{4.8cm}}$}}}
\put(16,10){El Ni\~no mode of 1997}
\put(16,7){(warm spot off Peru)}
\put(23,18){Year}
\put(4,40){\rotatebox{90}{Level}}
\put(96,68){$\psi^{(1,1)}_1$}
\put(59,74){mean temperature mode}
\put(96,50){$\psi^{(1,1)}_2$}
\put(96,34){$\psi^{(4,7)}_1$}
\put(95,30){$t\!\approx \!1997$}
\put(96,18){$\psi^{(4,9)}_1$}
\put(95,14){$t\! \approx \!1999$}
\end{overpic}
\end{center}
\vspace*{-.3in}
\caption{\label{fig:Fig4} Application of mrDMD on sea surface temperature data from 1990 to 2010. The left panel illustrates the process for a 4-level decomposition.  At each level, the {\em slowest} modes are extracted.  At a given level, the zero-mode component has a period $T=\infty$.  Illustrated are the (a) Level-1 mrDMD mode with period $T=\infty$, and (b) Level-1 mrDMD mode with period $T=52$ weeks.  These two modes are in the gray mode removal region of Fig.~\ref{fig:mrDMD}.  Further in the decomposition we can extract (c) Level-4 mrDMD mode of 1997 with period $T=\infty$ and (d) Level-4 mrDMD mode of 1999 with period $T=\infty$.   Mode (c) clearly shows
the El Ni\~no mode of interest that develops in the central and east-central equatorial Pacific.  The El Ni\~no mode was
not present in 1999 as is clear from mode (d).  Data source:  NOAA\_OI\_SST\_V2 data provided by the NOAA/OAR/ESRL PSD, Boulder, Colorado, USA, from their Web site at http://www.esrl.noaa.gov/psd/}
\end{figure}

The example of the previous subsection was contrived to demonstrate the ability of the mrDMD 
to separate the spatio-temporal modes of Fig.~\ref{fig:modes_temporal}.   In this example, we use
the mrDMD on a more realistic data set.  Specifically, we consider global surface temperature data
over the ocean.  The data is open source from the NOAA/OAR/ESRL PSD, Boulder, Colorado, USA.
The NOAA\_OI\_SST\_V2 data set considered can be downloaded at http://www.esrl.noaa.gov/psd/
The data spans a 20 year period from 1990 to 2010.  

Figure~\ref{fig:Fig4} shows the results of the mrDMD algorithm.  Specifically, a 4-level decomposition
is performed with the slow spatio-temporal modes pulled at each level as suggested in Fig.~\ref{fig:mrDMD}.
The zero mode in each window is the DC component, or period infinity $T=\infty$ mode.  At the
first level of the decomposition, two modes are extracted:  the zero mode $(T=\infty)$ depicted in Fig.~\ref{fig:Fig4}(a)
and a yearly cycle  $(T=52$ weeks) shown in Fig.~\ref{fig:Fig4}(b).  The yearly cycle is the {\em slowest} mode
extracted at level 1.  Note that the zero mode component of level 1, $\psi_1^{(1,1)}$, is just the average ocean temperature over the entire 20-year data set.

We can continue the mrDMD analysis through to the fourth-level.  At the fourth level, the 
approximate zero mode (period $T=\infty$) of the sampling window extracts physically
interesting results.  In particular, the data-driven method of the mrDMD discovers the 1997
El Ni{\~n}o mode generated from the well-known El Ni\~no, Souther Oscillation (ENSO).
Indeed, 1997 is known to have been a strong El Ni\~no year, as verified by its strong modal signature in the 4th level decomposition of the mrDMD.  In contrast, the same sampling
window shifted down to 1999 produces no El Ni\~no mode, which is in keeping with known
ocean patterns that year.   El Ni\~no is the warm phase of the ENSO cycle and is associated with a band of warm ocean water that develops in the central and east-central equatorial Pacific (between approximately the International Date Line and 120\degree W), including off the Pacific coast of South America.  The mrDMD mode clearly shows this
band of warm ocean water as a spatio-temporal mode in 1997 in Fig.~\ref{fig:Fig4}(c).  
ENSO refers to the cycle of warm and cold temperatures, as measured by sea surface temperature, SST, of the tropical central and eastern Pacific Ocean. El Ni\~no is accompanied by high air pressure in the western Pacific and low air pressure in the eastern Pacific. The cool phase of ENSO is called La Ni\~na with SST in the eastern Pacific below average and air pressures high in the eastern and low in western Pacific. The ENSO cycle, 
both El Ni\~no and La Ni\~na, causes global changes of both temperatures and rainfall.

These results could not have been produced with DMD
unless the correct sampling windows were chosen ahead of time, thus requiring a supervised learning step not required by the mrDMD.  Further, even in such a case,
the previous slow modes, such as those of level-1 demonstrated in Fig.~\ref{fig:Fig4}(a) and (b)
would pollute the data at the level of investigation.  Thus the mrDMD provides a principled, algorithmic
approach that is capable of data-driven discovery in complex systems such as ocean/atmospheric data.

\begin{figure}[t]
\begin{center}
\begin{overpic}[width=0.7\textwidth]{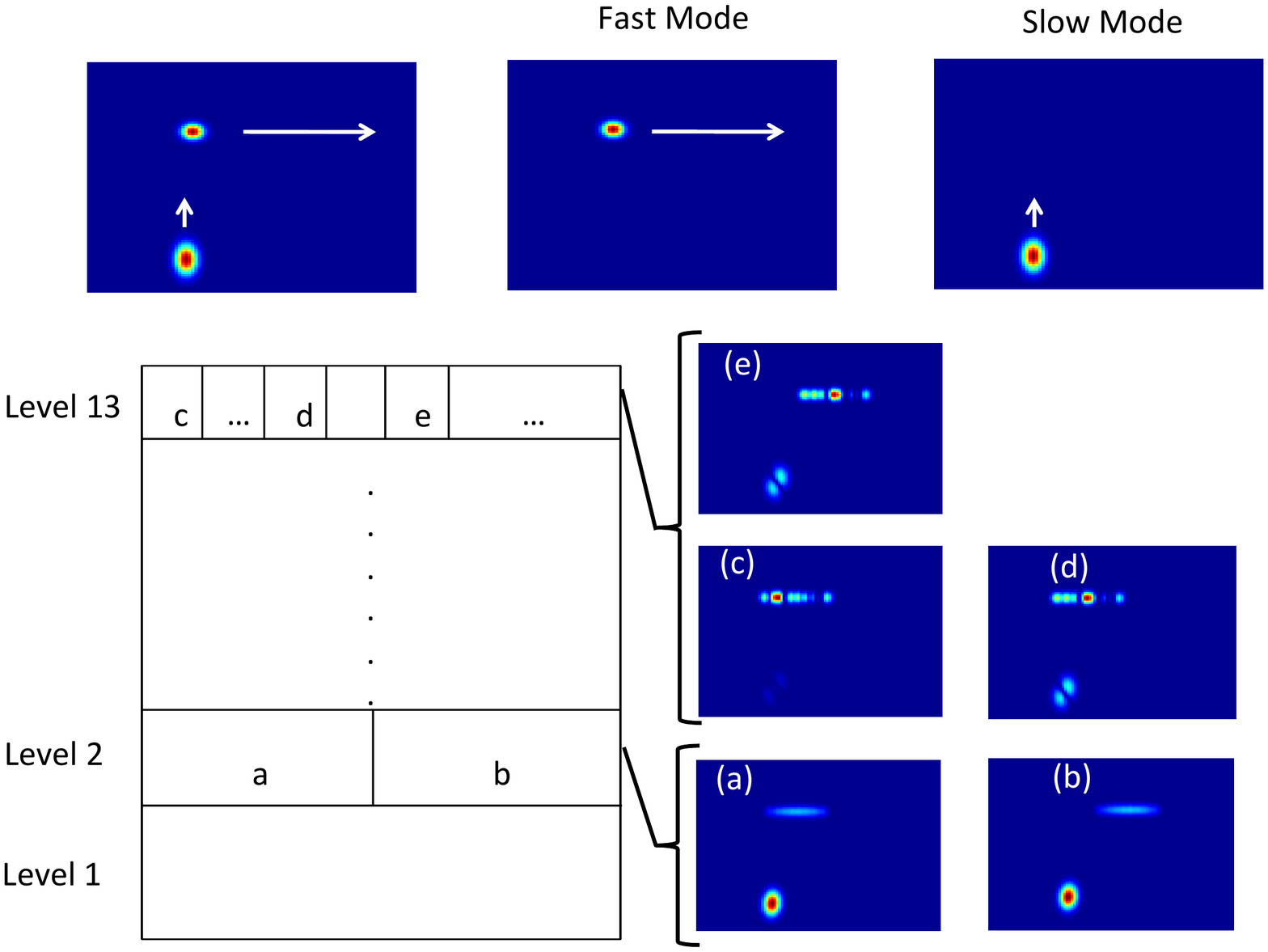}
\normalsize
\put(28.5,60){\color{white}{$10v$}}
\put(55,60){\color{white}{$10v$}}
\put(20,57){\color{white}{$v$}}
\put(75,57){\color{white}{$v$}}
\put(27,54){\color{white}{$\bar{\bf x} (x,y,t)$}}
\put(50,54){\color{white}{$\bar{\psi}_1 (x\!-\!10vt,\!y\!)$}}
\put(80,54){\color{white}{$\bar{\psi}_2 (x,\!y\!-\!vt\!)$}}
\put(40,59){$=$}
\put(67,59){$+$}
\put(64,13){\color{white}{$\psi^{(2,1)}_1$}}
\put(83,13){\color{white}{$\psi^{(2,2)}_1$}}
\put(64,27){\color{white}{$\psi^{(13,3)}_1$}}
\put(83,27){\color{white}{$\psi^{(13,5)}_1$}}
\put(64,40){\color{white}{$\psi^{(13,7)}_1$}}
\put(92,11){\rotatebox{90}{slow modes}}
\put(92,30){\rotatebox{90}{fast modes}}
\end{overpic}
\end{center}
\vspace*{-.5in}
\caption{\label{fig:Fig5} Application of mrDMD on moving objects separation. 
Two modes (top panel labeled ``Fast Mode" and ``Slow Mode") are combined in the
data snapshots (top right panel).   The ``Fast Mode" moves from left to right as indicated by
the arrow with a speed of $10v$ while the ``Slow Mode" moves from bottom to top at the
speed $v$.  We take $v=1/40$ without loss of generality.  
Thus the fast and slow mode speeds are approximately an order of magnitude
different.  In the mrDMD decomposition (left bottom panel), the ``Slow Mode" is extracted at level 2 
as represented by the panels (a) and (b).  The  ``Fast Mode" is extracted at level 13 as represented
at three representative panels (c)-(e).  Although there is some shadow (residual) of the slow mode on the fast mode and vice-versa, the mrDMD is remarkably effective in extracting the correct modes.  Moreover,
the level at which they are extracted can allow for a reconstruction of the velocity and direction.  
To our knowledge,
this is the best performance achieved to date with an SVD-based method with multiple time-scale objects.}
\end{figure}

\subsection{Translating and/or Rotating Structures}
The final application of the mrDMD is on an example that is
notoriously difficult for SVD-based methods to characterize, namely,  when
translational and/or rotational structures are present in the data snapshots, i.e. continuous or discrete 
invariances.   Indeed, such invariances
completely undermine our ability to compute low-rank embeddings of the
data as driven by correlated structures, or POD/PCA modes.   This has been
the Achilles heel of many SVD based methods, thus requiring in applications such
as PCA-based face recognition (eigenfaces)~\cite{eigenfaces}, well-cropped and centered faces for
reasonable performance, i.e. translation and rotation are removed in an expensive pre-processing
procedure.  

In a dynamical systems setting, a simple traveling wave will appear
to be a high-dimensional object in POD space despite the fact that it is only
constructed from two modes, one associated with translational invariance.
For dynamical cases exhibiting translation and/or rotation, Rowley and Marsden~\cite{Rowley20001} 
formulated one of the only mathematical strategies to date to extract
the low-dimensional embeddings.  In particular, they developed a template-matching technique 
to first remove the invariance before applying
the SVD decomposition.  Although effective, it is not suited for cases where a myriad of objects and timescales
are present in the data.

The mrDMD is well suited to handle invariances such as translation and rotation.  Consider
once again the case of a simple traveling wave.  The standard DMD decomposition 
applied to the traveling wave would result in a solution approximation requiring many 
DMD modes due to the slow fall off of the singular values (\ref{eq:kry3}) in step 1 of the DMD
algorithm.  Further, the eigenvalues $\omega_k$ in (\ref{eq:dmd_sol}) would also be bounded
away from the origin as there is no {\em background} mode for such translating data.
In the mrDMD architecture, the fact that the eigenvalues are bounded away from the origin
(and typically $O(1)$) in the initial snapshot window 
${\bf X}_M$ would simply allow the mrDMD method to ignore the traveling wave at the
first level of the mrDMD.  Once the domain is now divided into two for the next level of analysis, the
traveling wave is now effectively moving at half the speed in this new domain, i.e. the
eigenvalues have migrated towards the origin and the traveling wave is now re-evaluated.
The recursive procedure would eventually produce a sampling window where the traveling
wave looks sufficiently stationary and low-rank so as to be extracted in the multi-resolution analysis.
The level at which it is extracted also characterizes the speed of the traveling wave.  Specifically,
the higher the level in the decomposition where the traveling wave is extracted, the faster its
speed.

Figure~\ref{fig:Fig5} demonstrates the application of the mrDMD method on a simple
example in which there are two moving objects, one moving at a slow velocity and another
moving at a high velocity.  Specifically, the example constructed results from
the dynamical sequence
\begin{equation}
  \bar{\bf x} =  \bar{\psi}_1 (x-10vt,y) + \bar{\psi}_2 (x,y-vt)
\end{equation}
where the two modes (Gaussians of the form $\bar{\psi}_j=\exp(-\sigma(x-x_0)^2 - \sigma (y-y_0)^2)$ with
$\sigma=0.1$ and $(x_0,y_0)=(-18,20)$ and $(x_0,y_0)=(-20,-9)$ for the fast and slow modes respectively) used to construct the true solution are represented by $\bar{\psi}_j (x,y)$ for
$j=1,2$.  Note that the first mode is translating from left to right at speed $10v$ whereas the
second mode is translating from bottom to top at speed $v$.  Without loss of generality, $v$ can
be set to any value.  It is chosen to be $v=1/40$ for the domain and Gaussians considered.  In this case, the straightforward template matching procedure would fail due to
the two distinct time-scales of the objects.  As shown in the figure, the mrDMD is capable
of pulling out the two modes at level 2 (slow) and 13 (fast) of the analysis.  It is at these levels that the translating
objects look close to stationary ($|\omega_k|\ll 1$) in the mrDMD analysis.  The level 2 modes correspond
to the slow-moving object whereas the level 13 modes are associated with the fast object.
Although there is a residual in both the extracted slow and fast modes, it does a reasonable job
in extracting the fast and slow modes.  To our knowledge, this is
the only SVD-based method capable of such unsupervised performance.

The motivation for this example is quite clear when considering video processing
and surveillance.  In particular, for many applications in this arena, background subtraction
is only one step in the overall video assessment.  Identification of the foreground
objects is the next, and ultimately more important, task.  Consider the example of
a video of a street scene in which there is a pedestrian walking (slow translation) 
along with a vehicle driving down the road (fast translation).  Not only would we want
to separate the foreground from background, we would also want to separate the
pedestrian from the vehicle.  More precisely, we would want to 
make a separate video and identification of the slow and fast objects in the video.
The results of Fig.~\ref{fig:Fig5} show the mrDMD to be well-suited and effective in
this task.  Moreover, one can also envision augmenting the algorithm with a recursive
application of mrDMD to improve performance much like in the foreground/background
separation work on DMD~\cite{jake}.

\section{Discussion and Outlook}\label{sec:discussion}

Data-driven strategies for analyzing complex systems are of growing interest
in the mathematical sciences.   Indeed, methods capable of providing principled decompositions
of data arising in multi-scale, spatio-temporal systems are key enabling strategies
for many applications in the engineering, physical and biological sciences.
The multi-resolution dynamic mode decomposition advocated in this work capitalizes
on recent innovations in equation-free modeling via the dynamic mode decomposition.  It 
leverages these ideas by integrating them with concepts from wavelet theory and multi-resolution analysis.
By construction, the method provides a principled reconstruction of multi-resolution, spatio-temporal data sets.
The methods effectiveness is demonstrated on several example data sets, highlighting its
ability to extract critical information and enact data-driven discovery protocols.

The DMD method provides a spatio-temporal decomposition of data into
a set of dynamic modes that are derived from snap shots or measurements of a given system in time. 
The DMD method approximates the modes of the so-called Koopman operator. The Koopman operator is a linear, infinite-dimensional operator that represents nonlinear, infinite-dimensional dynamics without linearization, and is the adjoint of the Perron-Frobenius operator. The method can be viewed as computing, from the experimental data, the eigenvalues and eigenvectors (low-dimensional modes) of a linear model that approximates the underlying dynamics, even if the dynamics is nonlinear.  By interpreting the DMD eigenvalues as corresponding to prescribed time scale dynamics, one can
extract spatio-temporal structures recursively for shorter and shorter sampling windows.  Thus the slow-modes
are removed first and the data is filtered for analysis of its higher frequency content.  This recursive sampling
structure is demonstrated to be effective in allowing for a reconstruction of several example data sets.

One can envision a number of innovations to augment the proposed mrDMD strategy.  Many are
particularly attractive for applications across the engineering, physical and biological sciences.  
Indeed, the impact that DMD is having on complex fluid flows is already known and has already been mentioned in
the introduction.  Fields like neuroscience, which are rich in multi-scale, complex dynamics are also
ideal candidates for exploration using the mrDMD infrastructure as DMD has already had
recent demonstrated success in this arena~\cite{brunton2015}.
The outlook 
of these techniques is highlighted here:\\

\noindent
{\bf Compressive sampling:}  First is the ability to leverage tools from compressive sampling~\cite{cdmd} to facilitate the collection of considerably fewer measurements, resulting in the same multi-resolution dynamic mode decomposition, but with considerably fewer measurements. This reduction in the number of measurements may have a broad impact in situations where data acquisition is expensive and/or prohibitive. In particular, we envision these tools being used in particle image velocimetry (PIV) to reduce the data transfer requirements for each snapshot in time, increasing the maximum temporal sampling rate. Other applications include ocean and atmospheric monitoring, where individual sensors are expensive. Even if full-state measurements are available, the proposed method of compressed DMD will be computationally advantageous in situations where there is low-rank structure in the high-dimensional data.  

\noindent
{\bf DMD Control:}  A second important direction revolves around recent innovations
of the DMD with control~\cite{dmdc}, which is capable of disambiguating between the underlying dynamics and the effects of actuation, resulting in accurate input-output models. The method is data-driven in that it does not require knowledge of the underlying governing equations, only snapshots of state and actuation data from historical, experimental, or black-box simulations.  One can envision developing such input-output models at the various levels of spatio-temporal encoding in the data, i.e. multi-resolution input-output models can be assessed and constructed from a given system.  Both
these research directions, and many others, highlight the potential strength of the model in pushing forward
principled, equation-free strategies for the analysis of complex multi-scale systems.

\noindent
{\bf Machine Learning:}  Modern tools of statistical analysis and dimensionality-reduction have become the 
workhorses for the burgeoning field of machine learning (ML).   ML techniques aim to capitalize on
underlying low-dimensional patterns and clustering in data.  In the dynamical applications considered here, one might exploit these patterns, or DMD modes, by building libraries of low-rank dynamical modes, much like is done with POD modes~\cite{bright2013,brunton2014,epj2014}.  Such DMD libraries for different dynamical regimes partner nicely with compressive sensing strategies.  Additionally, Kernel based techniques, which are at the core
of support vector machines, for instance, have already found
successful application in the DMD architecture when considering more accurate, nonlinear dynamical 
reconstructions~\cite{williams2015}.  Maximum advantage should be
taken of such techniques when integrating the DMD architecture in applications.

\section*{Acknowledgment}
We are grateful for discussions with B. W. Brunton, J. Grosek and J. L. Proctor.  J. N. Kutz acknowledges support from the U.S. Air Force Office of Scientific Research (FA9550-09-0174).   (See http://youtu.be/E1dNE02LaCE
for a video summary)

\bibliographystyle{unsrt}
\bibliography{mrdmd.bib}
\end{document}